\documentclass[12pt,a4paper]{amsart}  
\usepackage[twosideshift=0cm,scale={0.73,0.8}]{geometry}
\newtheorem{theorem}{Theorem}[section]

\newtheorem{proposition}[theorem]{Proposition}   
\theoremstyle{remark}   
\newtheorem{remark}[theorem]{\sc Remark}   
\theoremstyle{definition}   
\newtheorem{definition}[theorem]{Definition}   
\theoremstyle{remark}   
\newtheorem{example}[theorem]{\sc Example}   
   
\newcommand{\cal}{\mathcal}    

\newcommand{\Aut}{\mathop{{\rm{Aut}}}\nolimits}   
\newcommand{\Reg}{\mathop{{\rm{Reg}}}\nolimits}   
\newcommand{\e}{\varepsilon}    
\newcommand{\ity}{{\infty}}   
\newcommand{\bR}{{\mathbb R}}   
\newcommand{\bC}{{\mathbb C}}   
\newcommand{\0}{{\bf{0}}}      
\newcommand{\cB}{{\cal B}}     
\newcommand{\cL}{{\cal L}}   
\newcommand{\cW}{{\cal W}}   
\begin{document}    
\title[Contact boundaries]{Contact boundaries of   
hypersurface singularities and of complex polynomials}     
\author{\sc Cl\'ement Caubel}  
\address{D{\'e}partement de Math{\'e}matiques, Universit{\'e} de  
Cergy-Pontoise, 95033 Cergy-Pontoise, France.}  
\email{caubel@math.u-cergy.fr}     
\author{\sc Mihai Tib\u ar}  
\address{Math\' ematiques, UMR-CNRS 8524, Universit\'e de Lille 1, \  
59655 Villeneuve d'Ascq, France.}     
\email{tibar@agat.univ-lille1.fr}  
\begin{abstract}   
We survey some recent results concerning the behaviour of the  
contact structure defined on the boundary of a complex isolated hypersurface singularity   
or on the boundary at infinity of a complex polynomial.  
\end{abstract}    
\maketitle   
\setcounter{section}{0}   
\section{Introduction}  
  
Let $f:(\bC^{n+1},\0)\to(\bC,0)$ be a germ of holomorphic function having  an isolated singular point at the origin. Its {\em boundary} $\cB_\0(f)$ is  
the intersection of the hypersurface $V(f):=f^{-1}(0)$ with a small sphere $S_\e$ centered  
at the origin, of radius $\e >0$. This is a closed oriented $(2n-1)$-dimensional  
smooth manifold, which does not depend on $\e\ll 1$ up to isotopy. If we consider  
its embedding in the sphere $S_\e$, that is, the  
{\em link} $\cL_\0(f):=(S_\e,S_\e\cap V(f))$,  
then it determines completely the topological type of $V(f)$  
(see \cite{Mi}). On the other hand, the natural CR-structure on $S_\e\cap V(f)$,   
defined by the maximal complex hyperplane distribution in its tangent bundle, determines  
its analytical type (see \cite{Scherk}). Thus we see that the boundary equiped with some  
additional structure can encode a lot of information on the singularity. 
 
If one considers a polynomial function $f:\bC^{n+1}\to\bC$ then we have to deal with more global objects. 
Replacing small spheres by large ones, the {\em boundary at infinity} of the general fiber of $f$ can be defined. 
Its embedding in the large sphere or its naturally defined CR-structure 
also provide potentially interesting invariants of $f$. 
 
We present here some recent results on an intermediate structure: the {\em contact structure} 
defined by the maximal complex hyperplane distribution. Section 2 collects the definitions and results in contact geometry that we use in the following.
 
In the local case, Varchenko \cite{Var} showed that this structure is an analytic invariant of 
 the singularity; we review this work in Section \ref{defi}. In the global case, we present the authors' construction \cite{CT} of the {\em contact boundary at infinity} attached to some polynomial. 
We therefore have well-defined contact boundaries attached to isolated singularities or to complex polynomials. 
 
These contact boundaries are analytic (respectively algebraic) invariants, but {\em a 
 priori} not topological invariants. This motivates the study of their variation in 
topologically trivial families of hypersurface singularities or of polynomials.  

Section 4 is devoted to recent results showing that the  
{\em formal homotopy class} (i.e. the most primitive contact invariant, see Example 
\ref{def-formal}) of these contact boundaries 
is constant in the following types of families, provided $n>2$:

\noindent

1.  topologically trivial families of isolated hypersurface singularities, \cite{Cau}.

\noindent

2. families of polynomials $f_s$ such that their general fibers are homotopy equivalent to a wedge of $n$-spheres and such that their 
$n^{\mbox{\footnotesize th}}$ betti numbers are constant in the family, \cite{CT}.
 
We end this survey by some remarks and questions in Section \ref{rems}. 
 
\section{A reminder on contact and almost contact structures}\label{remi} 
\subsection{Families of contact manifolds}   
For generalities on contact structures, see e.g.~\cite{Blair}, or \cite{El2}  
for a more up-to-date overview.  
  
\begin{definition}  
A {\em contact form} on a $(2n-1)$-dimensional manifold $M$ is a global 1-form  
$\alpha$ on $M$ satisfying the non-integrability condition $\alpha\wedge  
(d\alpha)^{n-1}>0$. A {\em contact structure} on $M$ is the hyperplane distribution  
defined by a contact form. The notions of {\em contactomorphism} or of {\em contact  
isotopy} are naturally defined.  
\end{definition}  
  
\begin{example}\label{pluricont} 
Any smooth level set of a strictly plurisubharmonic function on a complex manifold  
admits a natural contact structure, given by the maximal complex distribution in  
its tangent bundle  
(see e.g.~\cite{El1} for more details about contact structures and plurisubharmonic  
functions).  
\end{example}  
  
We now give the most general criterion for two contact structures to be isotopic.  
\begin{theorem}\label{gray}  
{\bf (J.W.~Gray \cite{Gray})} Let $\pi:M\to B$ be a smooth fiber bundle such that  
the fiber $M_*$ is a closed, oriented, odd-dimensional  
manifold. Suppose we are given a contact structure $\xi_b$ on each  
fiber $M_b$ which depends smoothly on the points $b$ of the base. Then  
the contact manifolds $(M_b,\xi_b)_{b\in B}$ are all contact isotopic.  
\end{theorem}  
  
This theorem will be used in Section \ref{defi} to show that the contact manifolds  
we associate to isolated singularities or complex polynomials are well defined.  
  
\subsection{Homotopies of almost contact manifolds}    
 
\begin{definition}  
Let $M$ be a smooth closed oriented odd-dimensional manifold. An {\em  
almost contact structure}\footnote{This definition  
differs slightly from the common one, see e.g.~\cite{Blair}.}  
on $M$ is a hyperplane distribution  
$\xi\subset TM$ endowed with a complex multiplication  
$J:\xi\to\xi$. An {\em almost contact manifold} is a manifold endowed  
with an almost contact structure. Finally, two almost contact manifolds  
$(M_1,\xi_1,J_1)$ and $(M_2,\xi_2,J_2)$ are {\em almost contact homotopic}  
if they are isotopic submanifolds of a manifold $W$ and if any isotopy between  
them carries $(\xi_1,J_1)$ on an almost contact structure on $M_2$  
which is homotopic to $(\xi_2,J_2)$ in the space of almost contact structures on  
$M_2$.  
\end{definition}    
  
\begin{example}\label{def-formal}  
Let $(M,\xi)$ be a contact manifold. Then any choice of a riemannian metric  
on $M$ and of a contact form defines a complex structure on $\xi$: the restriction  
of the 2-form $d\alpha$ on $\xi$ is non-degenerate, so it gives a complex structure  
on $\xi$ in presence of the metric. While this almost contact structure is not well  
defined, its homotopy class in  
the space of all almost contact structures on $M$ is: it is called the  
{\em formal homotopy class} of the contact structure $\xi$. This class  
is easily seen to be invariant up to contact  
isotopy, and provides the most primitive global invariant of contact structures (see  
\cite{El2}).  
\end{example}  
  
\begin{example}  
Any (real, cooriented) hypersurface $M$ in an almost complex manifold $(W,J)$ is  
naturally endowed with an almost contact structure: the complex  
hyperplane distribution $(\xi_M,J_M)$ is just the maximal complex subbundle  
$TM\cap J\,TM$ of the tangent bundle to the hypersurface equipped with the complex  
structure induced by $J$. Moreover, we have the following:  
\end{example}    
  
\begin{proposition}\label{prop-cle}  
Any two isotopic real hypersurfaces in an almost complex manifold  
are almost contact homotopic.   
\end{proposition}  
\begin{proof}[Sketch of proof]  
There are two equivalent ways to prove this. The first is to consider homotopy  
classes of almost contact structures on a $(2n-1)$-dimensional oriented manifold  
as group reductions of its tangent bundle from $SO_{2n-1}$ to $U_{n-1}$. Then  
in the context of this proposition, one builds a homotopy between the Gauss maps  
induced from the embeddings of the hypersurfaces in the almost complex manifold,  
which gives a homotopy of almost contact structures. This is the approach used in  
\cite{Cau}.  
  
The second way, used in \cite{CT}, is to identify almost contact structures with  
2-forms of maximal rank when a riemannian metric is fixed. One then chooses a  
hermitian metric on the ambiant almost complex manifold and constructs a homotopy  
between the associated 2-forms on the two hypersurfaces.  
\end{proof}  
  
This proposition means that if we define the  
{\em formal homotopy class} of a hypersurface in an almost complex  
manifold as the homotopy class of its almost contact structure, this  
formal homotopy class is an isotopy invariant. Of course, for a strictly  
pseudoconvex hypersurface (that is, a smooth  
level set of a strictly plurisubharmonic function) in a complex manifold, which  
is then also a contact manifold, the two notions of formal homotopy  
classes coincide.  
  
This proposition will be used in Section \ref{defo} to show to that the contact manifolds  
we associate to isolated singularities or complex polynomials stay in the same  
formal homotopy class in certain deformations.  
  
\section{Definition of the contact boundary in the local and global cases}\label{defi}  
\subsection{The local case: Varchenko's work}  
  
We summarize in this subsection the main results obtained by Varchenko in \cite{Var}.  
While it was known long before that the regular level sets of strictly  
plurisubharmonic functions on complex manifolds were natural examples of contact  
manifolds, it was the first time when this was used in the context of singularity  
theory.  
  
\begin{definition}  
Let $f:(\bC^{n+1},\0)\to(\bC,0)$ be a germ of holomorphic function with \0 as an  
isolated singular point. The {\em contact boundary} of $f$ is the boundary  
$\cB_\0(f)$ of $f$ equipped with the contact structure defined by the maximal  
complex distribution of its tangent bundle.  
\end{definition}  
  
The fact that this hyperplane distribution is indeed a contact structure is due  
to the strict plurisubharmonicity of the squared norm function $z\mapsto |z|^2$.  
  
\begin{proposition}{\rm(\cite{Var})}\label{var}
The contact boundary of an isolated hypersurface singularity does not depend on  
the radius of the small sphere defining it up to contact isotopy. Moreover, it  
neither depends on the choice of analytic coordinates at \0 in $\bC^{n+1}$.  
\end{proposition}  
  
\begin{proof}[Sketch of proof]  
Gray's theorem \ref{gray} is used two times. First, the squared norm function  
restricted to   
$V(f)$ is a proper submersion over $(0,\e_0)$ for $\e_0\ll 1$ sufficiently small:  
this defines a family of contact manifolds, hence a contact isotopy. Second, since  
the group of analytic changes of coordinates is path connected, and because of  
Whitney's (b)-property satisfied by $(V(f)\setminus \{\0\},\{\0\})$ (see \cite{Var}  
for details), one can construct a convenient family of contact manifolds and then  
conclude.  
\end{proof}  
  
This proposition shows that the contact boundary is a well-defined analytic  
invariant of the singularity defined by $f$.  
 
Having shown this invariance, Varchenko applies it to quasi-homogeneous  
singularities. He proves that in this case the contact boundary is {\em Sasakian}  
(a certain integrability condition, see e.g.~\cite{Blair}), which implies that its  
odd Betti numbers are even up to half its dimension.   
This gives a strong restriction for a singularity to be topologically equivalent  
to a quasi-homogeneous one. Varchenko considers as examples Brieskorn's cusp singularities:   
$$f_{p,q,r}:(x,y,z)\mapsto x^p+y^q+z^r+xyz,\qquad \frac{1}{p} +\frac{1}{q}
+\frac{1}{r}<1.$$  
Their boundaries satisfy $\dim H_1(\cB_\0(f_{p,q,r}))=1$, hence these singularities are
not topologically equivalent to quasi-homogeneous ones.  
  
  
\subsection{The global case}\label{ss:global} 
 Let us first take a complex algebraic set $V\subset\bC^{n+1}$ having at most 
isolated singularities. 
It is proved in \cite{CT} that, for large enough radius $R$,  
the intersection $S_R\cap V$ is a contact 
manifold, which is independent on $R$, up to contact isotopy. It turns out that  
instead of using spheres $S_R$, which are levels of the squared distance function,
one can use the levels of some {\em pseudo-convex rug function at infinity},
i.e.~a proper real polynomial map $\rho:\bC^{n+1}\to\bR_{\ge 0}$ which is strictly 
plurisubharmonic
(see Remark \ref{psh} for the relation with the local case).
Then the intersections $\rho^{-1}(R) \cap V$ do not depend on the
choice of $\rho$ and of $R\ge R_\rho\gg 1$, up to contact isotopy. These observations
yield a well-defined  {\em contact boundary of $V$} 
which we shall denote in the following by ${\cB}_\ity(V)$. 
 
Let now $f:\bC^{n+1}\to\bC$ be a complex polynomial function. It is 
well-known (see \cite{Th}, \cite{Br}) that there exists a  finite set 
$B_f\subset\bC$ (which we shall suppose minimal) such that the 
restriction $f_|:f^{-1}(\bC\setminus B_f)\to \bC\setminus B_f$ is a 
locally trivial  C$^\ity$ fibration. Any value $t\notin B_f$ is called 
{\em typical}, as is the corresponding fiber $f^{-1}(t)$. 
The atypical values (i.e.~those in $B_f$) are due not only to the 
critical points of the function $f$, but also to a certain bad 
asymptotic behaviour at infinity.  
 
It follows from the above discussion that any typical fiber  $f^{-1}(t)$ 
has a well-defined 
contact boundary ${\cB}_\ity(f^{-1}(t))$.   Moreover, we shall show how to define a {\em 
generic contact boundary}, independent on the fiber $f^{-1}(t)$, except for  
finitely many values of $t$.    
 
Let $\rho:\bC^{n+1}\to\bR_{\ge 0}$ be a pseudo-convex rug function at 
infinity. We say that the fiber $f^{-1}(t_0)$ is {\em 
$\rho$-regular-at-infinity} if there exists a (small enough) disk 
$D_\delta \subset \bC$ centered at $t_0$ and a (large enough) real 
$R_\delta \gg 0$ such that, for any $R\ge R_\delta$, the level 
$\rho^{-1}(R)$ is transversal to $f^{-1}(t)$ for all $t\in D_\delta$. 
 
It  follows from this definition that a $\rho$-regular-at-infinity 
fiber can have at most isolated singularities.  
By \cite{Ti1}, \cite[Prop.~2.6]{Ti2}, if   
 the fiber $f^{-1}(t_0)$ is non-singular and is 
$\rho$-regular-at-infinity for some $\rho$, then  $t_0$ is a typical value of $f$.\footnote{The converse of this statement is an open problem, to our 
knowledge (see \cite{ST}, \cite{Ti2}).}  
 
Let $\Reg_\ity f$ denote the set of values $t\in \bC$ such that the 
fiber  $f^{-1}(t)$ is $\rho$-regular-at-infinity, where $\rho$ is some pseudo-convex rug function at 
infinity.  
Then the precise meaning 
of the notion of ``generic contact boundary'' of a 
polynomial $f$ is given by the following. 
 
 
\begin{theorem}\rm (\cite{CT})\label{t:indep} \it 
Let $f:\bC^{n+1}\to\bC$ be a complex polynomial function. The contact boundary 
of a regular-at-infinity fiber $f^{-1}(t)$ does not 
depend, up to contact isotopy, on the choice of the value $t\in 
\Reg_\ity f$. 
\end{theorem}  
\begin{proof}[Sketch of proof]  
Let $d_0$ denote the squared distance function. We have proved in \cite[Cor.~2.12]{Ti2} that the set of values $t$ such that the fiber $f^{-1}(t)$ is not $d_0$-regular-at-infinity is a finite set. It then follows that the complement of $\Reg_\ity f$ in $\bC$ is a finite set (which of course  
contains $B_f$). 
 
We first reduce the problem to the case when $t_1,t_2 \in \Reg_\ity f$ 
are both $d_0$-regular-at-infinity values. Then use a path $\Gamma$ from $t_1$ to $t_2$  
such that $\Gamma$ consists of only  
$d_0$-regular-at-infinity values. For sufficiently  large $R$, 
$f_{| \Gamma} :f^{-1}(\Gamma)\cap S_R\to \Gamma$   
is a smooth fibration.  One applies Gray's Theorem \ref{gray} 
to show that the contact boundaries of the fibers $f^{-1}(t_1)$ and 
$f^{-1}(t_2)$ are contact isotopic. 
\end{proof}  
 
 
The above theorem defines an invariant which we call {\em contact boundary of $f$} and which we denote by 
$\cB_\ity f$. Let us remark that it does not depend on changes of coordinates in the 
automorphism group $\Aut\bC^{n+1}$ since the class of pseudo-convex rug 
functions at infinity is right-invariant by the action of $\Aut\bC^{n+1}$.  
 
  
\section{Evolution of the contact boundary  
in numerically constant deformations}\label{defo}  
\subsection{The local case: invariance of the formal homotopy class in  
$\mu$-constant deformations}  
  
Let $f:(\bC^{n+1},\0)\to(\bC,0)$ be a holomorphic germ defining an isolated singularity at  
\0. In \cite{Cau}, the first named author adresses the following question: how does  
the contact boundary vary in topologically trivial deformations $f$? A partial answer was 
obtained (see also Remark \ref{eliash} on Eliashberg's question): the most primitive global invariant of the contact boundary is preserved in  
such deformations. More precisely, we have the following:

\begin{theorem}{\rm (\cite{Cau})}\label{defo-loca}  
Let $n>2$ and let $(f_s)_{s\in[0,1]}:(\bC^{n+1},\0)\to(\bC,0)$ be a smooth family  
of holomorphic  
function germs having an isolated singular point at \0. Suppose that all these germs  
are topologically equivalent. Then the contact boundaries $\cB_\0 f_0$ and  
$\cB_\0 f_1$ are almost contact homotopic.  
\end{theorem}

\begin{proof}[Sketch of proof]  
We rely on the proof in \cite{LR} of the L\^e-Ramanujam theorem relating topological  
triviality and   
$\mu$-constancy of deformations. For small enough, well-chosen positive real parameters  $s,t,\e_0, \e_s$, one can use Gray's Theorem \ref{gray} to show  
the following contact isotopies:  
$$\cB_\0 f_0 = f_0^{-1}(0)\cap S_{\e_0} \stackrel{\rm cont}{\simeq}  
f_0^{-1}(t)\cap S_{\e_0}  
\stackrel{\rm cont}{\simeq} f_s^{-1}(t)\cap S_{\e_0}:=M_0$$  
$$\cB_\0 f_s = f_s^{-1}(0)\cap S_{\e_s} \stackrel{\rm cont}{\simeq}  
f_s^{-1}(t)\cap S_{\e_s}=: M_s.$$   
Denote by  
$W:=f_s^{-1}(t)\cap\left(B_{\e_0}\setminus  
\stackrel{\circ}{B_{\e_s}}\right)$ the piece of the fiber  
$f_s^{-1}(t)$ between the two spheres of radius $\e_0$ and  
$\e_s$. The same reasoning as in \cite{LR} shows that, under our  
hypotheses, $W$ is a $h$-cobordism between $M_0$ and $M_s$. Since $n>2$, the two  
hypersurfaces  
$M_0\simeq\cB_\0 f_0$ and $M_s\simeq\cB_\0 f_s$  
are isotopic in the complex manifold $W$. We then use Proposition \ref{prop-cle}  
to conclude.  
\end{proof}

\subsection{The global case: invariance of the formal homotopy class in  
$\gamma$-constant V-deformations}  
We shall call {\em family of polynomials} a family $f_s : \bC^{n+1}  
\to \bC$ depending smoothly on the parameter $s\in [0, \e]$, for some  
real $\e > 0$.  
We have used in the local case the equivalence between the topological triviality  
and the $\mu$-constancy, which allowed us to conclude by using the $h$-cobordism theorem. 
In case of polynomial functions, even if we suppose that $f_s$ has isolated singularities, the total Milnor number is by far not enough to control the topology of $f_s$. A class of polynomials which turns out to be a good 
candidate for the analogue of the class of germs of isolated singularities
is given in the following. 
  
\begin{definition}\label{d:v}  
We say that a family of polynomials is a {\em V-deformation of $f =  
f_0$} if, for all $s\in [0, \e]$, the following condition is  
satisfied:  
$$ \mbox{(V)\quad the typical fiber of $f_s$ is homotopy equivalent to a bouquet of spheres } \bigvee_{\gamma_s} S^n.$$   
If moreover the number $\gamma_s$ of spheres is independent on $s$,  
then we say that we have a {\em $\gamma$-constant V-deformation}.  
\end{definition} 
  
  The class of $\gamma$-constant V-deformations has been considered in  
\cite{Ti1, Ti2} for proving a L\^e-Ramanujam type result concerning the constancy of the monodromy at infinity. It is easy to see that if a V-deformation is topologically trivial then it is $\gamma$-constant. The converse is not true, see \cite{Ti2}, \cite{CT}.

\begin{example}\label{e:v}   
  The polynomials with isolated  
$\cW$-singularities at infinity in the sense of Siersma-Tib\u ar \cite{ST}  are examples of polynomials satisfying condition (V). This class includes  
the polynomials such that all their fibers are $\rho$-regular-at-infinity (for some $\rho$) and all polynomials of 2 complex  
variables with irreducible fibers.   
\end{example}   
  
\begin{theorem}\label{defo-glob}{\rm (\cite{CT})}  
Let $n>2$ and let $(f_s)_{s\in[0,\e]}:\bC^{n+1}\to\bC$ be a   
$\gamma$-constant V-deformation of polynomial functions.  
Then the contact boundaries $\cB_\ity f_0$ and  
$\cB_\ity f_1$ are almost contact homotopic. In particular,  
this is the case in topologically trivial V-deformations.  
\end{theorem}  
  
\begin{proof}[Sketch of proof]  
It has been observed in \cite{Ti1, Ti2}   
that the arguments used in the local  
case by L\^e-Ramanujam (\cite{LR}) can be applied  
in case of polynomials satisfying the condition  
(V) when replacing the small spheres by large ones.  
In our context, we paraphrase the proof of Theorem \ref{defo-loca}  
like this: for a smooth $d_0$-regular at infinity fiber $f_0^{-1}(t)$  
and for a good choice of parameters satisfying  
$1\gg R_0^{-1}\gg s\gg \delta\gg R_s^{-1}>0$, we get the  
following contact isotopies:  
$$\cB_\ity f_0 \stackrel{\rm cont}{\simeq} f_0^{-1}(t)\cap S_{R_0}  
\stackrel{\rm cont}{\simeq} f_s^{-1}(t)\cap S_{R_0} \stackrel{\rm  
cont}{\simeq} f_s^{-1}(t+\delta)\cap S_{R_0},$$  
$$\cB_\ity f_s\stackrel{\rm cont}{\simeq}f_s^{-1}(t+\delta)\cap S_{R_s}.$$  
The cobordism  
$W:=f_s^{-1}(t+\delta)\cap\left(B_{R_s}\setminus  
\stackrel{\circ}{B_{R_0}}\right)$  
between the two contact boundaries is a product $h$-cobordism under  
our hypotheses, which allows us to conclude using Proposition \ref{prop-cle}.  
Notice that we were bound to pass from $t$ to $t+\delta$ in the base, since   
the $d_0$-regularity at infinity of $f_s^{-1}(t)$ is not guaranteed.  
\end{proof}  
  
\section{Remarks and questions}\label{rems} 
\begin{remark}\label{psh}
For simplicity, we have defined in the local case the boundary of an analytic
function germ by means of the squared norm function. But this choice may seem
too restrictive: for instance, in \cite{Lo} the boundary in defined by any
{\em rug function}, that is, any proper
real analytic function $\rho:(\bC^{n+1},\0)\to\bR_{\geq 0}$ such that
$\rho^{-1}(0)=\0$. This doesn't change the isotopy type of the boundary.
Since we want to get a well-defined contact structure on the boundary, the
good generalization in our context is the notion of {\em pseudo-convex rug function},
that is, a strictly plurisubharmonic rug function\footnote{
This was used by Ehlers, Neumann and Scherk in \cite{ENS} for surface singularities}.
This doesn't change the contact isotopy type of the contact boundary,
which gives another proof of Proposition \ref{var}: the class of pseudo-convex rug
functions is right invariant by the action of the group $\Aut(\bC^{n+1},\0)$
of analytic changes of coordinates. This is the approach we have used in the global
case.
\end{remark}
  
\begin{remark}\label{eliash}  
In \cite{El1}, Eliashberg raises the following ``$J$-convex $h$-cobordism problem'':  
suppose we have a strictly plurisubharmonic function $f$ on a product  
$W=M\times[0,1]$, such that the two boundary components $M\times\{i\}$, $i=0,1$,  
are level sets of $f$. Is it possible to find another strictly plurisubharmonic  
function $g$ which coincides with $f$ on the boundary and without critical points?  
  
Were this true, the proof of Theorem \ref{defo-loca} would improve to yielding that  
 the  
contact boundary is invariant in topologically trivial deformations. 
In the global case too (Theorem \ref{defo-glob})
this would yield the invariance of the contact boundaries in  $\gamma$-constant (and, in particular,
topologically trivial) V-deformations.  
   
\end{remark}    
  
\begin{remark}  
The results of Section \ref{defo} were proved using the $h$-cobordism theorem,  
which excludes the surface case $n=2$. For normal surface germs,  
thanks to the plumbed structure of the boundary given by its topological
type (see \cite{N}), one can show that the formal homotopy class of the contact
boundary is
in fact a topological invariant. This will appear elsewhere \cite{C2}.  
\end{remark}

\begin{remark}  
Following the preceding Remark, one may ask in the general case if the formal
homotopy class
of the contact boundary is a topological invariant, which is {\em a priori}
a stonger condition than being an invariant of topologically  
trivial deformations. 
The answer is only known in the particular case when this boundary is  
diffeomorphic to the standard sphere $S^{2n-1}$. This is {\em yes}, in the
local case, by Morita's work \cite{Morita}. Morita gives formulae  
expressing the formal homotopy class (which is then an element of  
$\pi_{2n-2}(SO_{2n-1}/U_{n-1})$) in terms of the Milnor number $\mu$ and of the  
signature of the Milnor fiber, which are both topological invariants.  
\end{remark}

\end{document}